\documentclass[12pt]{amsart}
\usepackage{amssymb}
\usepackage[margin=0.8in]{geometry}

\usepackage[hyphenbreaks]{breakurl}
\usepackage[hyphens]{url}
\newcommand{\nc}{\newcommand}
\nc{\set}[2]{\{\, #1 : #2\,\}}
\nc{\st}{\op{star}}
\nc{\tPa}[6]{\bibitem{#1} {#2}, \emph{#3}, {#4}, to appear (#5 pages).}
\nc{\sPa}[5]{\bibitem{#1} {#2}, \emph{#3}, {#4}, submitted.}
\nc{\FinSeqs}[1]{{{#1}^{<\alephes}}}
\nc{\rest}[1]{\mid_{#1}}
\nc{\Bgp}{{\Z^\N}}
\nc{\Arh}{Arhangel'ski\u{\i}}
\nc{\grbl}{{\mbox{\textit{\tiny gp}}}}
\long\def\forget#1\forgotten{}
\nc{\issuenumber}{39}
\nc{\issuemonth}{July}
\nc{\issueyear}{'16}

\nc{\Ga}{\Gamma}
\nc{\sub}{\subseteq}
\nc{\R}{\mathbb{R}}
\nc{\Cp}{\mathrm{C}_\mathrm{p}}
\nc{\Op}{\mathrm{O}}
\nc{\Fin}{{[\N]^{<\aleph_0}}}
\nc{\alephes}{{\aleph_0}}
\nc{\ed}{

\forget

\section{Unsolved problems from earlier issues}

\begin{issue}Is $\binom{\Omega}{\Gamma}=\binom{\Omega}{\Tau}$?\end{issue}
\begin{issue}Is $\ufin(\Op,\Omega)=\sfin(\Gamma,\Omega)$?And if not, does $\ufin(\Op,\Gamma)$ imply
$\sfin(\Gamma,\Omega)$?\end{issue}
\stepcounter{issue}
\begin{issue}Does $\sone(\Omega,\Tau)$ imply $\ufin(\Gamma,\Gamma)$?\end{issue}
\stepcounter{issue}
\begin{issue}Is there, in ZFC, an uncountable set satisfying $\sfin(\cB,\cB)$?\end{issue}
\stepcounter{issue}
\begin{issue}Does $X \nin \NON(\cM)$ and $Y\nin\mathsf{D}$ imply that $X\cup Y\nin \COF(\cM)$?\end{issue}
\begin{issue}[CH]Is $\split(\Lambda,\Lambda)$ preserved under finite unions?\end{issue}
\begin{issue}Is $\cov(\cM)=\fo$? (See the definition of $\fo$ in that issue.)\end{issue}
\stepcounter{issue}
\begin{issue}Could there be a Baire metric space $M$ of weight $\aleph_1$ and a partition
$\mathcal{U}$ of $M$ into $\aleph_1$ meager sets where for each ${\mathcal U}'\subset\mathcal U$,
$\bigcup {\mathcal U}'$ has the Baire property in $M$?\end{issue}
\stepcounter{issue} 
\begin{issue}Is there, in ZFC, a set of reals $X$ of cardinality $\fd$ such that all
finite powers of $X$ have Menger's property $\sfin(\Op,\Op)$?\end{issue}
\stepcounter{issue}
\begin{issue}[MA]Is there a set $X\sbst\bbR$, of cardinality continuum, satisfying $\sone(\BO,\BG)$?\end{issue}
\begin{issue}[CH]Is there a totally imperfect $X$ satisfying $\ufin(\Op,\Gamma)$
that can be mapped continuously onto $\Cantor$?\end{issue}
\begin{issue}[CH]Is there a Hurewicz $X$ such that $X^2$ is Menger but not Hurewicz?\end{issue}
\begin{issue}Does the Pytkeev property of $C_p(X)$ imply that $X$ has Menger's property?\end{issue}
\begin{issue}Does every hereditarily Hurewicz space satisfy $\sone(\BG,\BG)$?\end{issue}
\begin{issue}[CH]Is there a Rothberger-bounded $G\le\Bgp$ such that $G^2$ is not Menger-bounded?\end{issue}
\begin{issue}Let $\cW$ be the van der Waerden ideal. Are $\cW$-ultrafilters closed under products?\end{issue}
\begin{issue}Is the $\delta$-property equivalent to the $\gamma$-property $\binom{\Omega}{\Gamma}$?\end{issue}
\stepcounter{issue}\stepcounter{issue}

\forgotten

\general\end{document}}

\nc{\fbx}[1]{\fbox{$#1$}}\nc{\nop}{$\times$}\nc{\fbn}{\!\!\fbox{\!\nop\!}\!\!}
\nc{\yup}{\checkmark}\nc{\fby}{\!\!\fbox{\!\yup\!}\!\!}\nc{\mbq}{\mb{?}}
\nc{\mb}[1]{{\mbox{\textbf{#1}}}}\nc{\smb}[1]{{\!\!\mb{#1}\!\!}}\nc{\x}{\times}
\nc{\<}{\left <}\nc{\Cantor}{{\{0,1\}^\N}}\nc{\oo}{\infty}
\nc{\NR}{{\bbR^\N}}\nc{\Iff}{\Leftrightarrow}\nc{\mypar}[1]{\par\medskip\noindent\textbf{#1.}}
\nc{\roth}{{[\N]^{\alephes}}}\nc{\sr}[2]{{\txt{$#1$\\$#2$}}}
\nc{\gpbl}{{\mbox{\textit{\tiny gp}}}}\nc{\fx}{\mathfrak{x}}\nc{\fb}{\mathfrak{b}}
\nc{\fg}{\mathfrak{g}}\nc{\fc}{\mathfrak{c}}\nc{\fd}{\mathfrak{d}}
\nc{\fp}{\mathfrak{p}}\nc{\fs}{\mathfrak{s}}\nc{\ADD}{{\mathsf   {ADD}}}
\nc{\COV}{{\mathsf   {COV}}}\nc{\NON}{{\mathsf   {NON}}}\nc{\COF}{{\mathsf   {COF}}}
\nc{\sseq}[1]{\{#1 : n\in\N\}}\nc{\Impl}{\Rightarrow}
\nc{\upannouncement}[1]{[\S\ref{#1} above]}\nc{\dnannouncement}[1]{[\S\ref{#1} below]}
\nc{\E}{\exists}\nc{\cI}{\mathcal{I}}\nc{\cN}{\mathcal{N}}\nc{\cP}{\mathcal{P}}
\nc{\cA}{\mathcal{A}}\nc{\cM}{\mathcal{M}}\nc{\Null}{\mathcal{N}}
\nc{\op}{\operatorname}\nc{\cov}{\mathsf{cov}}\nc{\add}{\mathsf{add}}
\nc{\cof}{\mathsf{cof}}\nc{\cf}{\mathsf{cf}}\nc{\non}{\mathsf{non}}\nc{\spst}{\supseteq}
\nc{\CH}{the Continuum Hypothesis}\nc{\bbR}{\mathbb{R}}\nc{\Q}{\mathbb{Q}}
\nc{\EdNote}[1]{\par\medskip\noindent\textbf{#1.}}\nc{\fo}{\mathfrak{od}}
\nc{\cl}[1]{\overline{#1}}\nc{\impl}{\rightarrow}\nc{\arrays}{{{\{0,1\}}^{\N\x\N}}}
\nc{\w}{\omega}\nc{\ft}{\mathfrak{t}}\nc{\h}{\mathfrak{h}}\nc{\Cite}[2]{{\cite[#1]{#2}}}
\renewcommand{\split}{\mathsf{Split}}\nc{\bq}{\begin{quote}}\nc{\eq}{\end{quote}}
\nc{\cK}{\mathcal{K}}\nc{\cB}{\mathcal{B}}\nc{\BG}{\cB_\Gamma}
\nc{\BL}{\cB_\Lambda}\nc{\BT}{\cB_\Tau}\nc{\BTstar}{\cB_{\Tau^*}}\nc{\BO}{\cB_\Omega}
\nc{\CG}{C_\Gamma}\nc{\CL}{C_\Lambda}\nc{\CT}{C_\Tau}\nc{\CTstar}{C_{\Tau^*}}
\nc{\CO}{C_\Omega}\nc{\sone}{\mathsf{S}_1}\nc{\sfin}{\mathsf{S}_\mathrm{fin}}
\nc{\Sc}{\mathsf{S}_c}\nc{\ufin}{\mathsf{U}_\mathrm{fin}}\nc{\gone}{\mathsf{G}_1} \nc{\gfin}{\mathsf{G}_\mathrm{fin}}\nc{\seq}[1]{\{#1\}_{n\in\N}}\nc{\Un}{\bigcup}
\nc{\nin}{\not\in}\nc{\cF}{\mathcal{F}}\nc{\cG}{\mathcal{G}}\nc{\cU}{\mathcal{U}}
\nc{\cV}{\mathcal{V}}\nc{\cW}{\mathcal{W}}\nc{\fU}{\mathfrak{U}}\nc{\fu}{\mathfrak{u}}
\nc{\fV}{\mathfrak{V}}\nc{\fW}{\mathfrak{W}}\nc{\psin}{pseudo-intersection}
\nc{\NN}{{\N^\N}}\nc{\N}{\mathbb{N}}\nc{\bbN}{\mathbb{N}}\nc{\Z}{\mathbb{Z}}
\nc{\as}{\subseteq^*}\nc{\sm}{\setminus}\nc{\sbst}{\subseteq}
\nc{\by}[2]{\par\hfill\emph{#1}, #2}\nc{\nby}[1]{\par\hfill\emph{#1}}\nc{\Tau}{\mathrm{T}}
\nc{\CE}{\textsc{CE}}
\newtheorem{thm}{Theorem}[section]\nc{\bthm}{\begin{thm}} \nc{\ethm}{\end{thm}}
\newtheorem{prop}[thm]{Proposition}\nc{\bprp}{\begin{prop}} \nc{\eprp}{\end{prop}}
\newtheorem{fact}[thm]{Fact}\nc{\bfct}{\begin{fact}} \nc{\efct}{\end{fact}}
\newtheorem{prob}[thm]{Problem}\nc{\bprb}{\begin{prob}} \nc{\eprb}{\end{prob}}
\newtheorem{lem}[thm]{Lemma}\nc{\blem}{\begin{lem}} \nc{\elem}{\end{lem}}
\newtheorem{claim}[thm]{Claim}\nc{\bclm}{\begin{claim}} \nc{\eclm}{\end{claim}}
\newtheorem{cor}[thm]{Corollary}\nc{\bcor}{\begin{cor}} \nc{\ecor}{\end{cor}}
\newtheorem{conj}[thm]{Conjecture}\nc{\bcnj}{\begin{conj}} \nc{\ecnj}{\end{conj}}
\theoremstyle{definition}\newtheorem{defn}[thm]{Definition}\nc{\bdfn}{\begin{defn}} \nc{\edfn}{\end{defn}}
\theoremstyle{remark}\newtheorem{rem}[thm]{Remark}\nc{\brem}{\begin{rem}} \nc{\erem}{\end{rem}}
\newtheorem{cnv}[thm]{Convention}\nc{\bcnv}{\begin{cnv}} \nc{\ecnv}{\end{cnv}}
\newtheorem{exam}[thm]{Example}\nc{\bexm}{\begin{exam}} \nc{\eexm}{\end{exam}}
\newtheorem{issue}{Issue}\nc{\bpf}{\begin{proof}} \nc{\epf}{\end{proof}}
\nc{\be}{\begin{enumerate}}\nc{\ee}{\end{enumerate}}\nc{\bi}{\begin{itemize}}
\nc{\ei}{\end{itemize}}\nc{\itm}{\item}
\nc{\general}{\small\vfill\par\noindent\hrulefill\par
\noindent\textbf{Previous issues.} 
\url{http://front.math.ucdavis.edu/search?\&t=\%22SPM+Bulletin\%22}
\\[0.1cm]
\textbf{Contributions and free subscription.} Email \url{tsaban@math.biu.ac.il}.
}

\nc{\link}[1]{\par\hfill{\url{#1}}}
\nc{\fillin}{{\Huge To be completed}}
\nc{\arXivl}[4]{\subsection{#2}{#4}\par\hfill{\arx{#1}}\par\hfill\emph{#3}}
\nc{\arXiv}[3]{\subsection{#2} {#3}, \arx{#1}\par\hfill}
\nc{\DOIpaper}[5]{\subsection{#2}{#4}\par\hfill{\url{http://dx.doi.org/#1}}\par\hfill\emph{#3}}
\nc{\AMSPaper}[5]{\subsection{#3}{#5}\par\hfill{\url{#1}}\par\hfill\emph{#4}\par\hfill{#2}}
\nc{\nAMSPaper}[4]{\subsection{#2} {#3}, {#4}, \url{#1}}
\nc{\AMS}[3]{\subsection{#1} {#2}, \url{#3}}
\nc{\SPMBul}{\textbf{$\mathcal{SPM}$ Bulletin}}
\nc{\BulEnd}{\par\bigskip\noindent
Boaz Tsaban\\
\emph{E-mail}: tsaban@math.biu.ac.il\\
\emph{URL}: http://www.cs.biu.ac.il/\~{}tsaban}

\nc{\arx}[1]{\url{http://arxiv.org/abs/#1}}

\nc{\probissue}{\emph{Problem of the issue}}

\title[$\mathcal{SPM}$ Bulletin \textbf{\issuenumber} (\issuemonth{} \issueyear)]{%
$\mathcal{SPM}$ Bulletin\\[0.5cm]
Issue \issuenumber{} (\issuemonth{} \issueyear{})}

\begin{document}
\maketitle


\section{Editor's note}

\medskip

With the approaching TOPOSYM'16 (\url{http://www.toposym.cz/programme.php}), 
it is a pleasure to 
see how the topic of selection principles gains increasing
attention and becomes a standard part of topology and set theory.
At least eight of the 28 speakers, and a good number 
of the contributed lecture speakers, made substantial contributions to
this topic in their career. 
For some of these, SPs constitute the main topic
of research in the last few years. 
This is in accordance with the continuous progress on the topic, 
some of which mentioned below.

With best regards,

\by{Boaz Tsaban}{tsaban@math.biu.ac.il}

\hfill \texttt{http://www.cs.biu.ac.il/\~{}tsaban}

\section{Long announcements}

\arXivl{1512.00788}
{A characterization of barrelledness of $C_p(X)$}
{S. Gabriyelyan}
{We prove that, for a Tychonoff space $X$, the space $C_p(X)$ is barrelled if
and only if it is a Mackey group.}

\arXivl{1512.07515}
{Countable tightness and $\mathfrak G$-bases on Free topological groups}
{Fucai Lin, Alex Ravsky, Jing Zhang}
{Given a Tychonoff space $X$, let $F(X)$ and $A(X)$ be respectively the free
topological group and the free Abelian topological group over $X$ in the sense
of Markov. In this paper, we discuss two topological properties in $F(X)$ or
$A(X)$, namely the countable tightness and $\mathfrak G$-base. We provide some
characterizations of the countable tightness and $\mathfrak G$-base of $F(X)$
and $A(X)$ for various special classes of spaces $X$. Furthermore, we also
study the countable tightness and $\mathfrak G$-base of some level of $F(X)$.
Some open problems in \cite{GKL} are partially answered.}

\arXivl{1601.01010}
{Ideal equal Baire classes}
{Adam Kwela and Marcin Staniszewski}
{ For any Borel ideal we characterize ideal equal Baire system generated by the
families of continuous and quasi-continuous functions, i.e., the families of
ideal equal limits of sequences of continuous and quasi-continuous functions.}

\arXivl{1601.03163}
{On questions which are connected with Talagrand problem}
{Volodymyr V. Mykhaylyuk}
{We prove the following results.
\be
\item If $X$ is a $\alpha$-favourable space, $Y$ is a regular space, in which
every separable closed set is compact, and $f:X\times Y\to\mathbb R$ is a
separately continuous everywhere jointly discontinuous function, then there
exists a subspace $Y_0\subseteq Y$ which is homeomorphic to $\beta\mathbb N$.
\item There exist a $\alpha$-favourable space $X$, a dense in $\beta\mathbb
N\setminus\mathbb N$ countably compact space $Y$ and a separately continuous
everywhere jointly discontinuous function $f:X\times Y\to\mathbb R$.
\ee  
  Besides, it was obtained some conditions equivalent to the fact that the
space $C_p(\beta\mathbb N\setminus\mathbb N,\{0,1\})$ of all continuous
functions $x:\beta\mathbb N\setminus\mathbb N\to\{0,1\}$ with the topology of
point-wise convergence is a Baire space.
}

\arXivl{1601.05884}
{Namioka spaces and strongly Baire spaces}
{Volodymyr V. Mykhaylyuk}
{A notion of strongly Baire space is introduced. Its definition is a
transfinite development of some equivalent reformulation of the Baire space
definition. It is shown that every strongly Baire space is a Namioka space and
every $\beta-\sigma$-unfavorable space is a strongly Baire space.}

\arXivl{1602.06227}
{Comparing Fr\'echet-Urysohn filters with two pre-orders}
{S. Garcia--Ferreira, J. E. Rivera--G\'omez}
{A filter $\cF$ on $\w$ is called Fr\'echet-Urysohn if the space with only one
non-isolated point $\w \cup \{\cF\}$ is a Fr\'echet-Urysohn space, where the
neighborhoods of the non-isolated point are determined by the elements of $\cF$.
In this paper, we distinguish some Fr\'echet-Urysohn filters by using two
pre-orderings of filters: One is the Rudin-Keisler pre-order and the other one
was introduced by Todor\v{c}evi\'c-Uzc\'ategui in \cite{tu05}. In this paper,
we construct an RK-chain of size $\fc^+$ which is RK-above of avery
FU-filter. Also, we show that there is an infinite RK-antichain of
FU-filters.}

\arXivl{1603.09715}
{Selective game versions of countable tightness with bounded finite
  selections}
{Leandro F. Aurichi, Angelo Bella, Rodrigo R. Dias}
{For a topological space $X$ and a point $x \in X$, consider the following
game -- related to the property of $X$ being countably tight at $x$. In each
inning $n\in\omega$, the first player chooses a set $A_n$ that clusters at $x$,
and then the second player picks a point $a_n\in A_n$; the second player is the
winner if and only if $x\in\overline{\{a_n:n\in\omega\}}$.
  In this work, we study variations of this game in which the second player is
allowed to choose finitely many points per inning rather than one, but in which
the number of points they are allowed to choose in each inning has been fixed
in advance. Surprisingly, if the number of points allowed per inning is the
same throughout the play, then all of the games obtained in this fashion are
distinct. We also show that a new game is obtained if the number of points the
second player is allowed to pick increases at each inning.}

\arXivl{1604.02473}{Gruff Ultrafilters}
{David Fern\'andez-Bret\'on and Michael Hru\v{s}\'ak}
{We investigate the question of whether $\mathbb Q$ carries an ultrafilter
generated by perfect sets (such ultrafilters were called gruff ultrafilters by
van Douwen). We prove that one can (consistently) obtain an affirmative answer
to this question in three different ways: by assuming a certain parametrized
diamond principle, from the cardinal invariant equality $\mathfrak d=\mathfrak
c$, and in the Random real model.}

\arXivl{1604.02555}{Free locally convex spaces with a small base}
{Saak Gabriyelyan and Jerzy Kakol}
{The paper studies the free locally convex space $L(X)$ over a Tychonoff space
$X$. Since for infinite $X$ the space $L(X)$ is never metrizable (even not
Fr\'echet-Urysohn), a possible applicable generalized metric property for
$L(X)$ is welcome. We propose a concept (essentially weaker than
first-countability) which is known under the name a $\mathfrak{G}$-base. A
space $X$ has a {\em $\mathfrak{G}$-base} if for every $x\in X$ there is a base
$\{ U_\alpha : \alpha\in\mathbb{N}^\mathbb{N}\}$ of neighborhoods at $x$ such
that $U_\beta \subseteq U_\alpha$ whenever $\alpha\leq\beta$ for all
$\alpha,\beta\in\mathbb{N}^\mathbb{N}$, where
$\alpha=(\alpha(n))_{n\in\mathbb{N}}\leq \beta=(\beta(n))_{n\in\mathbb{N}}$ if
$\alpha(n)\leq\beta(n)$ for all $n\in\mathbb{N}$. We show that if $X$ is an
Ascoli $\sigma$-compact space, then $L(X)$ has a $\mathfrak{G}$-base if and
only if $X$ admits an Ascoli uniformity $\mathcal{U}$ with a
$\mathfrak{G}$-base. We prove that if $X$ is a $\sigma$-compact Ascoli space of
$\mathbb{N}^\mathbb{N}$-uniformly compact type, then $L(X)$ has a
$\mathfrak{G}$-base. As an application we show: (1) if $X$ is a metrizable
space, then $L(X)$ has a $\mathfrak{G}$-base if and only if $X$ is
$\sigma$-compact, and (2) if $X$ is a countable Ascoli space, then $L(X)$ has a
$\mathfrak{G}$-base if and only if $X$ has a $\mathfrak{G}$-base.}

\arXivl{1604.04609}
{On separability of the functional space with the open-point and
  bi-point-open topologies, II}
{Alexander V. Osipov}
{In this paper we continue to study the property of separability of functional
space $C(X)$ with the open-point and bi-point-open topologies.}

\arXivl{1604.05116}
{On sequential separability of functional spaces}
{Alexander V. Osipov, Evgenii G. Pytkeev}
{In this paper, we give necessary and sufficient conditions for the space
$B_1(X)$ of first Baire class functions on a Tychonoff space $X$, with pointwise
topology, to be (strongly) sequentially separable.}

\arXivl{1604.08868}
{No interesting sequential groups}
{Alexander Shibakov}
{We prove that it is consistent with ZFC that no sequential topological groups
of intermediate sequential orders exist. This shows that the answer to a 1981
question of P.~Nyikos is independent of the standard axioms of set theory. The
model constructed also provides consistent answers to several questions of
D.~Shakhmatov, S.~Todor\v{c}evi\'c and Uzc\'ategui. In particular, we show that
it is consistent with ZFC that every countably compact sequential group is
Fr\'echet-Urysohn.}

\arXiv{1604.08872}
{On sequential analytic groups}
{Alexander Shibakov}
{We answer a question of S.~Todor\v{c}evi\'c and C.~Uzc\'ategui from
\cite{TU1} by showing that the only possible sequential orders of sequential
analytic groups are $1$ and $\omega_1$. Other results on the structure of
sequential analytic spaces and their relation to other classes of spaces are
given as well. In particular, we provide a full topological classification of
sequential analytic groups by showing that all such groups are either
metrizable or $k_\omega$-spaces, which, together with a result by Zelenyuk,
implies that there are exactly $\omega_1$ non homeomorphic analytic sequential
group topologies.}

\arXivl{1604.08874}
{On large sequential groups}
{Alexander Shibakov}
{We construct, using $\diamondsuit$, an example of a sequential group $G$ such
that the only countable sequential subgroups of $G$ are closed and discrete,
and the only quotients of $G$ that have a countable pseudocharacter are
countable and Fr\'echet. We also show how to construct such a $G$ with several
additional properties (such as make $G^2$ sequential, and arrange for every
sequential subgroup of $G$ to be closed and contain a nonmetrizable compact
subspace, etc.).
  Several results about $k_\omega$ sequential groups are proved. In particular,
we show that each such group is either locally compact and metrizable or
contains a closed copy of the sequential fan. It is also proved that a dense
proper subgroup of a non Fr\'echet $k_\omega$ sequential group is not
sequential extending a similar observation of T.~Banakh about countable
$k_\omega$ groups.}

\arXivl{1606.01013}
{The Ascoli property for function spaces}
{Saak Gabriyelyan, Jan Greb\'ik, Jerzy Kakol, Lyubomyr Zdomskyy}
{The paper deals with Ascoli spaces $C_p(X)$ and $C_k(X)$ over Tychonoff
spaces $X$. The class of Ascoli spaces $X$, i.e. spaces $X$ for which any
compact subset $K$ of $C_k(X)$ is evenly continuous, essentially includes the
class of $k_{\mathbb R}$-spaces. First we prove that if $C_p(X)$ is Ascoli,
then it is $\kappa$-Fr\'echet-Urysohn. If $X$ is cosmic, then $C_p(X)$ is
Ascoli iff it is $\kappa$-Fr'echet-Urysohn. This leads to the following
extension of a result of Morishita: If for a \v{C}ech-complete space $X$ the
space $C_p(X)$ is Ascoli, then $X$ is scattered. If $X$ is scattered and
stratifiable, then $C_p(X)$ is an Ascoli space. Consequently: (a) If $X$ is a
complete metrizable space, then $C_p(X)$ is Ascoli iff $X$ is scattered. (b) If
$X$ is a \v{C}ech-complete Lindel\"of space, then $C_p(X)$ is Ascoli iff $X$ is
scattered iff $C_p(X)$ is Fr\'echet-Urysohn. Moreover, we prove that for a
paracompact space $X$ of point-countable type the following conditions are
equivalent: (i) $X$ is locally compact. (ii) $C_k(X)$ is a $k_{\mathbb
	R}$-space. (iii) $C_k(X)$ is an Ascoli space. The Asoli spaces $C_k(X,[0,1])$
are also studied.}

\arXivl{1603.03361}
{Products of Menger spaces: a combinatorial approach}
{Piotr Szewczak, Boaz Tsaban}
{We construct Menger subsets of the real line whose product is not Menger in
	the plane. In contrast to earlier constructions, our approach is purely
	combinatorial, and the set theoretic hypotheses used are either milder than or
	incompatible with earlier ones. On the other hand, we establish productive
	properties for versions of Menger's property parameterized by filters and
	semifilters. In particular, the hypothesis $\mathfrak{b}=\mathfrak{d}$ implies
	that every productively Menger set of real numbers is productively Hurewicz,
	and each ultrafilter version of Menger's property is strictly between Menger's
	and Hurewicz's classic properties. We include a number of open problems
	emerging from this study.}

\arXivl{1607.01687}
{Products of Menger spaces, II: general spaces}
{Piotr Szewczak, Boaz Tsaban}
{We study products of general topological spaces with Menger's covering
property, and its refinements based on filters and semifilters.
Among other results, we prove that, assuming the Continuum Hypothesis, every
productively Lindel\"of space is productively Menger, and every productively
Menger space is productively Hurewicz. None of these implications is
reversible.}

\arXivl{1607.03599}
{Topological spaces with a local $\omega^\omega$-base have the strong
Pytkeev$^*$ property}
{Taras Banakh}
{Modifying the known definition of a Pytkeev network, we introduce a notion of
Pytkeev$^*$ network and prove that a topological space has a countable Pytkeev
network at a point $x\in X$ if and only if $X$ is countably tight at $x$ and
has a countable Pykeev$^*$ network at $x$. We define a topological space $X$ to
have the strong Pytkeev$^*$ property if $X$ has a Pytkeev$^*$ network at each
point. Our main theorem says that a topological space $X$ has a countable
Pytkeev$^*$ network at point $x\in X$ if $X$ has a local $\omega^\omega$-base
at $x$ (i.e., a neighborhood base $(U_\alpha)_{\alpha\in\omega^\omega}$ at $x$
such that $U_\beta\subset U_\alpha$ for all $\alpha\le\beta$ in
$\omega^\omega$). Consequently, each countably tight space with a local
$\omega^\omega$-base has the strong Pytkeev property.}


\section{Short announcements}\label{RA}

\arXiv{1512.02458}
{Luzin $\pi$-bases and operation foliage hybrid}
{Mikhail Patrakeev}

\arXiv{1601.02343}
{Separately continuous functions on products and its dependence on
  $\aleph$ coordinates}
{Volodymyr V. Mykhaylyuk}

\arXiv{1601.04081}
{Reversible filters}
{Alan Dow, Rodrigo Hern\'andez-Guti\'errez}

\arXiv{1603.05366}
{The Bi-Compact-Open Topology on C($X$)}
{A. Jindal, R. A. McCoy, S. Kundu}

\arXiv{1603.07459}
{Completeness and compactness properties in metric spaces, topological
  groups and function spaces}
{Alejandro Dorantes-Aldama, Dmitri Shakhmatov}

\arXiv{1604.04005}
{Free topological vector spaces}
{Saak S. Gabriyelyan and Sidney A. Morris}

\arXiv{1604.06178}
{Continuous extension of functions from countable sets}
{V. Mykhaylyuk}

\AMS{The left side of Cicho\'n's diagram}
{Martin Goldstern; Diego Alejandro Mejia; Saharon Shelah}
{http://www.ams.org/journal-getitem?pii=S0002-9939-2016-13161-4}

\arXiv{1605.00853}
{A primitive associated to the Cantor-Bendixson derivative on the real
  line}
{Borys \'Alvarez-Samaniego and Andr\'es Merino}

\arXiv{1605.01024}
{On Borel semifilters}
{Andrea Medini}

\arXiv{1605.04087}
{Every filter is homeomorphic to its square}
{Andrea Medini, Lyubomyr Zdomskyy}

\arXiv{1606.01967}
{$\mathfrak G$-bases in free (locally convex) topological vector spaces}
{Taras Banakh and Arkady Leiderman}

\arXiv{1607.01491}
{Completeness Properties of the open-point and bi-point-open topologies
on C$(X)$}
{Anubha Jindal, R. A. McCoy, S. Kundu and Varun Jindal}

\arXiv{1607.00517}
{On $\sigma$-countably tight spaces}
{Istv\'an Juh\'asz and Jan van Mill}

\ed